\documentclass{ifacconf}

\usepackage{graphicx}      
\usepackage{natbib}        
\usepackage{amsmath}
\usepackage{amssymb}
\usepackage{mathdots}
\usepackage[ruled,vlined,algo2e]{algorithm2e}
\usepackage[all,cmtip]{xy}
\usepackage{tikz}
\usetikzlibrary{automata,positioning,shapes,arrows,shadows}
\begin{document}
\begin{frontmatter}

\pdfminorversion=4

\title{Computing Semilinear Sparse Models for Approximately Eventually Periodic Signals}

\author[Fredy]{Fredy Vides} 

\address[Fredy]{
Scientific Computing Innovation Center,
Universidad Nacional Aut\'onoma de Honduras, Tegucigalpa, Honduras, (e-mail: fredy.vides@unah.edu.hn)
}

\begin{abstract}
Some elements of the theory and algorithmics corresponding to the computation of semilinear sparse models for discrete-time signals are presented. In this study, we will focus on approximately eventually periodic discrete-time signals, that is, signals that can exhibit an aperiodic behavior for an initial amount of time, and then become approximately periodic afterwards. The semilinear models considered in this study are obtained by combining sparse representation methods, linear autoregressive models and GRU neural network models, initially fitting each block model independently using some reference data corresponding to some signal under consideration, and then fitting some mixing parameters that are used to obtain a signal model consisting of a linear combination of the previously fitted blocks using the aforementioned reference data, computing sparse representations of some of the matrix parameters of the resulting model along the process. Some prototypical computational implementations are presented as well.
\end{abstract}

\begin{keyword}
Autoregressive models, neural-network models, parameter identification, least-squares approximation, time-series analysis.
\end{keyword}

\end{frontmatter}

\section{Introduction}
In this document, some elements of the theory and algorithmics corresponding to the computation of semilinear sparse models for discrete-time signals are presented. The study reported in this document is focused on approximately eventually periodic discrete-time signals, that is, signals that can exhibit an aperiodic behavior for an initial amount of time, and then become approximately periodic afterwards.

The main contribution of the work reported in this document is the application of a {\em colaborative scheme} involving sparse representation methods, linear autoregressive models and GRU neural network models, where each block model is first fitted independently using some reference data corresponding to some given signal. Subsequently, some mixing parameters are fitted, in order to obtain a signal model consisting of a linear combination of the previously fitted blocks, using the aforementioned reference data to fit the mixing coefficients. Along the process, some of the matrices of parameters of the resulting model are fitted using sparse representation methods. Some theoretical aspects of the aforementioned process are described in \S\ref{sec:semilinear-approximation}. As a byproduct of the work reported in this document, a toolset of Python programs for semilinear sparse signal model computation based on the ideas presented in \S\ref{sec:semilinear-approximation} and \S\ref{sec:algorithm} has been developed, and is available in \cite{FVides_SPAAR}.

The applications of the sparse signal model identification technology developed as part of the work reported in this document, range from numerical simulation for predictive maintenance of industrial equipement and structures, to geological data analysis. 

A prototypical algorithm for the computation of sparse semilinear autoregressors based on the ideas presented in \S\ref{sec:semilinear-approximation}, is presented in \S\ref{sec:algorithm}. Some illustrative computational implementations of the prototypical algorithm presented in \S\ref{sec:algorithm} are documented in \S\ref{sec:experiments}.

\section{Preliminaries and Notation}
Given $\delta>0$, let us consider the function defined by the expression 
\begin{align*}
H_\delta(x)=\left\{
\begin{array}{ll}
1, & x>\delta\\
0,& x\leq \delta
\end{array}
\right..
\end{align*}
Given a matrix $A\in \mathbb{C}^{m\times n}$ with singular values denoted by the expressions $s_j(A)$ for $j=1,\ldots,\min\{m,n\}$. We will write $\mathrm{rk}_\delta(A)$ to denote the number
\begin{align*}
\mathrm{rk}_\delta(A)=\sum_{j=1}^{\min\{m,n\}}H_\delta(s_j(A)).
\end{align*}

Given a time series $\Sigma=\{x_{t}\}_{t\geq 1}\subset \mathbb{C}$, a positive integer $L$ and any $t\geq L$, we will write $\mathbf{x}_L(t)$ to denote the vector
\[
\mathbf{x}_L(t)=\begin{bmatrix}
x_{t-L+1} & x_{t-L+2} & \cdots & x_{t-1}  &  x_{t}
\end{bmatrix}^\top \in \mathbb{C}^{L}.
\]

Given an ordered sample $\Sigma_{N}=\{x_{t}\}_{t=1}^{N}\subset \Sigma$ from a time series $\Sigma=\{x_t\}_{t\geq 1}$, we will write $\mathcal{H}_L(\Sigma_N)$ to denote the Hankel type trajectory matrix corresponding to $\Sigma_{N}$, that is determined by the following expression.
\begin{align*}
\mathcal{H}_L(\Sigma_N)&=\begin{bmatrix}
x_1 & x_2 & x_3 & \cdots & x_{N-L+1}\\
x_2 & x_3 & x_4 & \cdots & x_{N-L+2}\\
\vdots & \vdots & \vdots & \iddots & \vdots\\
x_{L} & x_{L+1} & x_{L+2} & \cdots & x_{N} 
\end{bmatrix}
\end{align*}

We will write $I_n$ to denote de identity matrix in $\mathbb{C}^{n\times n}$, and we will write $\hat{e}_{j,n}$ to denote the matrices in $\mathbb{C}^{n\times 1}$ representing the canonical basis of $\mathbb{C}^{n}$ (each $\hat{e}_{j,n}$ equals the $j$-column of $I_n$).

We will write $\mathbf{S}^1$ to denote the set $\{z\in \mathbb{C}:|z|=1\}$. Given any matrix $X\in \mathbb{C}^{m\times n}$, we will write $X^\ast$ to denote the conjugate transpose $\overline{X}^\top\in \mathbb{C}^{n\times m}$ of $X$. A matrix $P\in \mathbb{C}^{n\times n}$ will be called an orthogonal projector whenever $P^2=P=P^\ast$. Given any matrix $A\in \mathbb{C}^{n\times n}$, we will write $\sigma(A)$ to denote the spectrum of $A$, that is, the set of eigenvalues of $A$.

\section{Semilinear Modeling of Approximately Eventually Periodic Signals} \label{sec:semilinear-approximation}
A discrete-time signal represented by a times series $\Sigma=\{x_{t}\}_{t\geq 1}$ is said to be approximately eventually periodic ({\bf AEP}) if it can be aperiodic for an initial amount of time, and then becomes approximately periodic afterwards. In other words, there are $\varepsilon>0$ and two integers $T,S>0$ such that $|x_{t+kT}-x_{t}|\leq \varepsilon$ for each $t\geq S$ and each integer $k\geq 0$. The integer $T$ will be called an approximate period of $\Sigma$.

\begin{rem}\label{rem:AEP-rem}
Based on the notion of approximately eventually periodic signal considered on this study, it can be seen that given an AEP signal $\Sigma=\{x_t\}_{t\geq 1}$, there is a positive integer $S$ such that the tail $\{x_{t}\}_{t\geq S}$ of $\Sigma$ is approximately periodic.
\end{rem}

\subsection{Semilinear Sparse Autoregressors}
\label{autoregressors}
Given an AEP signal $\Sigma=\{x_t\}_{t\geq 1}\subset \mathbb{C}$ and a lag value $L>0$. For the study reported in this document we will consider semilinear signal models of the form:
\begin{align}
x_{t+1}&=\mathcal{S}(\mathbf{x}_L(t))\nonumber\\
&=\mathcal{L}(\mathbf{x}_L(t))+\mathcal{G}(\mathbf{x}_L(t))+\mathcal{E}(\mathbf{x}_L(t)), t\geq L,
\label{eq:ARGRU_model}
\end{align}
where $\mathcal{L}(\mathbf{x}_L(t))$ denotes a linear operation defined by the expression
\begin{equation}
\mathcal{L}(\mathbf{x}_L(t))=c_{1}x_{t}+c_{2}x_{t-1}+c_{3}x_{t-2}+\cdots+c_{L}x_{t-L+1},
\label{eq:AR-part}
\end{equation}
the term $\mathcal{G}(\mathbf{x}_L(t))$ represents a linear combination of neural networks based on gated recurrent unit ({\bf GRU}) layers, whose structure is described by the block diagram
\begin{equation}
\tikzstyle{block} = [draw, fill=white, rectangle, 
    minimum height=3em, minimum width=6em]
\tikzstyle{sum} = [draw, fill=white, circle, node distance=1cm]
\tikzstyle{input} = [coordinate]
\tikzstyle{output} = [coordinate]
\tikzstyle{pinstyle} = [pin edge={to-,thin,black}]
\begin{tikzpicture}[auto, node distance=2cm,>=latex']

    \node [input, name=input] {};
    \node (u0) [coordinate, xshift = -3.3cm , yshift = 0cm] {};
    \node (u1) [coordinate, yshift = -0.3cm] {};
    
    \node [block, right of=input, node distance=2cm] (K) {$\mathbf{A}$};
    \node [block, left of=input, node distance=1.2cm] (L) {$\mathbf{GRU}$};

    \draw [->] (input) -- node[name=u] {$\mathbf{h}(t)$} (K);   
    \draw [->] (u0) -- node[name=u] {$\mathbf{x}_L(t)$} (L);

    \node [output, right of=K] (output) {};
    \draw [->] (K) -- node [name=y] {$\hat{x}_{t+1}$}(output); 
    \node [near end] {$ $} (input); 
\end{tikzpicture}
\label{eq:Adapted_GRU_diagram}
\end{equation}
and $\mathcal{E}(\mathbf{x}_L(t))$ represents some suitable error term. For some given integer $m>0$, each GRU cell $\mathbf{G}_j$ in the GRU block in \eqref{eq:Adapted_GRU_diagram} is described for each $j=1,\ldots,m$ by the following equations:
\begin{align}
r_j(t) &= \sigma\left(\hat{e}_{j,m}^\top \left(W_{ir}\mathbf{x}_L(t)+W_{hr}\mathbf{h}(t-1)+b_{r}\right)\right)\nonumber\\
z_j(t) &= \sigma\left(\hat{e}_{j,m}^\top \left(W_{iz}\mathbf{x}_L(t)+W_{hz}\mathbf{h}(t-1)+b_{z}\right)\right)\nonumber\\
n_j(t) &=\tanh(\hat{e}_{j,m}^\top \left(W_{in}\mathbf{x}_L(t)+b_{n}\right)\nonumber\\
&+r_j(t)\hat{e}_{j,m}^\top \left(W_{hn}\mathbf{h}(t-1)\right))\nonumber\\
h_j(t)&=(1-z_j(t))n_j(t)+z_j(t) h_j(t-1)
\label{eq:GRU-cell-description}
\end{align}
with $\mathbf{h}(t)=\begin{bmatrix}
h_{1}(t) & \cdots & h_{m}(t)
\end{bmatrix}^\top$, and where $\sigma$ denotes the sigmoid function. The configuration considered in \eqref{eq:Adapted_GRU_diagram} is based on a GRU layer in order to prevent vanishing gradients, by taking advantage of the GRU structure presented in \cite{cho-etal-2014-learning}. Although LSTM networks can be used to prevent vanishing gradients as well, for the study reported in this document GRU networks were chosen instead of LSTM networks, as they have fewer trainable parameters. The affine layer $\mathbf{A}$ of the neural network described in \eqref{eq:Adapted_GRU_diagram} is determined by the expression
\[
\mathbf{A}(\mathbf{h}(t))=\mathbf{w}_A^\top \mathbf{h}(t)+b_A.
\]

In order to compute models of the form \eqref{eq:ARGRU_model}, we can combine sparse autoregressive models of the form \eqref{eq:AR-part} that can be computed using the methods presented in \cite{DBLP:journals/corr/abs-2105-07522}, with GRU neural network models of the form \eqref{eq:Adapted_GRU_diagram} that can be computed using the computational tools provided as part of TensorFlow, Keras and PyTorch, that are described as part of \cite{chollet2015keras} and \cite{NEURIPS2019_9015}.

An approximate representation 
\[
\tilde{\mathcal{L}}(\mathbf{x}_L(t))=\tilde{c}_{1}x_{t}+\tilde{c}_{2}x_{t-1}+\tilde{c}_{3}x_{t-2}+\cdots+\tilde{c}_{L}x_{t-L+1},
\]
of the linear component of \eqref{eq:ARGRU_model} such that
\begin{align*}
\mathcal{L}(\mathbf{x}_L(t))\approx \tilde{\mathcal{L}}(\mathbf{x}_L(t)), t\geq L
\end{align*}
can be computed using some sample $\Sigma_N=\{x_t\}_{t=1}^N$ and a corresponding subsample $\Sigma_{0}=\{x_{t}\}_{t=1}^{N-1}\subset \Sigma_N$ for some suitable $N>L$, by {\em approximately solving} the matrix equation
\begin{align}\label{eq:sp_autoregressor_eq}
\mathcal{H}_L(\Sigma_{0})^\top \begin{bmatrix}
c_{L} \\ c_{L-1} \\ \vdots \\ c_{2} \\ c_{1}
\end{bmatrix}=
\begin{bmatrix}
x_{L+1}\\
x_{L+2}\\
\vdots\\
x_{N-1}\\
x_{N}
\end{bmatrix},
\end{align}
using sparse least-squares approximation methods. 

Schematically, the semilinear autoregressors considered in this study can be described by a block diagram of the form,
\begin{equation}
\tikzstyle{sensor}=[draw, text width=3em, 
    text centered, minimum height=2.5em]
\tikzstyle{ann} = [above, text width=5em, text centered]
\tikzstyle{wa} = [sensor, text width=2em, 
    minimum height=4em, rounded corners]
\tikzstyle{sc} = [sensor, text width=13em, fill=red!20, 
    minimum height=10em, rounded corners, drop shadow]
\def\blockdist{2.3}        
\begin{tikzpicture}
    \node (wa) [wa]  {$\mathfrak{M}$};
    \path (wa.west)+(-1.5,1.5) node (asr1) [sensor] {$\mathfrak{L}$};
    \path (wa.west)+(-1.5,0.5) node (asr2)[sensor] {$\mathfrak{G}_1$};
    \path (wa.west)+(-1.5,-1.0) node (dots)[ann] {$\vdots$}; 
    \path (wa.west)+(-1.5,-2.0) node (asr3)[sensor] {$\mathfrak{G}_N$};       
    \path (wa.east)+(1.1,0) node (vote) {$x_{t+1}$};
    \path [draw, ->] (asr1.east) -- node [above] {} 
        (wa.160) ;
    \path [draw, ->] (asr2.east) -- node [above] {} 
        (wa.180);
    \path [draw, ->] (asr3.east) -- node [above] {} 
        (wa.200);
    \path [draw, ->] (wa.east) -- node [above] {} 
        (vote.west);
    \path [draw, ->] (asr1.west)+(-1,0) node [left] {$\mathbf{x}_L(t)$} -- node {} (asr1.west);        
    \draw [->] (asr1.west)+(-0.5,0) node {} |- node {} (asr2.west);
    \draw [->] (asr2.west)+(-0.5,0) node {} |- node {} (asr3.west);
\end{tikzpicture}
\label{eq:Model_diagram}
\end{equation}
where the block $\mathfrak{L}$ is represented by \eqref{eq:AR-part}, each block $\mathfrak{G}_j$ is represented by \eqref{eq:Adapted_GRU_diagram}, and where the block $\mathfrak{M}$ is a {\em mixing} block defined by the expression
\begin{equation*}
\mathfrak{M}(y_1(t),\ldots,y_{N+1}(t))=\sum_{j=1}^{N+1} w_{j}y_j(t),
\end{equation*}
for some coefficients $w_j$ to be determined and some given $N$, with $y_1(t)=\mathfrak{L}(\mathbf{x}_L(t))$ and $y_{k+1}(t)=\mathfrak{G}_k(\mathbf{x}_L(t))$ for each $k=1,\ldots,N$ and each $t\geq L$.

The details of the computation of the neural network blocks of model \eqref{eq:ARGRU_model} will be omitted for brevity, for details on the theory and computation of the GRU neural network models considered for this study the reader is kindly referred to \cite{cho-etal-2014-learning}, \cite{chollet2015keras}, \cite{NEURIPS2019_9015} and \cite{FVides_SPAAR}.

Several interesting papers have been written on the subject of hybrid time series models that combine ARIMA and ANN models, two important references on this matter are \cite{ZHANG2003159} and \cite{KHANDELWAL2015173}. Besides using sparse AR models instead of ARIMA models, another important distinctive aspect of the modeling approach reported in this document, is that instead of using the recurrent neural network components of \eqref{eq:Model_diagram} represented by $\mathcal{G}$ in \eqref{eq:ARGRU_model} to approximate the residuals $r_t=x_{t+1}-\mathcal{L}(\mathbf{x}_L(t))$. Using some suitable training subsets $\Sigma_{I},\Sigma_M$ of a given data sample $\Sigma_{N}$ from an arbitrary AEP signal $\Sigma=\{x_t\}_{t\geq 1}$ under consideration, first the parameters of the blocks $\mathfrak{L}$, $\mathfrak{G}_1$, $\cdots$, $\mathfrak{G}_N$ of \eqref{eq:Model_diagram} are fitted independently using $\Sigma_I$, and then the coefficients of the block $\mathfrak{M}$ of \eqref{eq:Model_diagram} are fitted using $\Sigma_M$ and some corresponding predicted values generated with $\mathfrak{L}$, $\mathfrak{G}_1$, $\cdots$, $\mathfrak{G}_N$. Computing sparse representations of some of the matrix parameters of the resulting model along the process.

\subsection{An Operator Theoretic Approach to the Computation \\
of Linear Components of Semilinear Sparse Autoregressors}

Given an AEP signal $\Sigma=\{x_{t}\}_{t\geq 1}$ whose behavior can be approximately described by a model of the form \eqref{eq:AR-part}, that can be computed by approximately solving an equation of the form \eqref{eq:sp_autoregressor_eq} for some suitable integers $N,L>0$ with $N>L$, and given some sample $\Sigma_0=\{x_t\}_{t=1}^{N-1}$, if we consider any sample $\tilde{\Sigma}_{0}=\{\tilde{x}_t\}_{t=1}^{N-1}\subset \Sigma$, such that for some positive integer $S$ the states in $\tilde{\Sigma}_{0}$ satisfy the conditions $\tilde{x}_t=x_{t+S}$, for each $t=1,\ldots,N-1$. We will have that the matrix 
\begin{align}
C_L=\begin{bmatrix}
0 & 1 & 0 & \cdots & \cdots & 0\\
0 & \ddots & \ddots & \ddots &  & \vdots\\
\vdots & \ddots & \ddots & \ddots & \ddots & \vdots\\
\vdots &  & \ddots & \ddots & \ddots & 0\\
0 & \cdots & \cdots & 0 & 0 & 1\\
c_L & c_{L-1} & \cdots & \cdots & c_{2} & c_{1}
\end{bmatrix}\in \mathbb{C}^{L\times L}
\label{eq:C_matrix}
\end{align}
will {\em approximately} satisfy the condition
\begin{align}\label{eq:sp_autoregressor_matrix_eq}
\mathcal{H}_L(\Sigma_{0})^\top \left(C_L^S\right)^\top=
\mathcal{H}_L(\tilde{\Sigma}_{0})^\top.
\end{align}

Using matrices of the form \eqref{eq:C_matrix} one can express linear models of the form \eqref{eq:AR-part} as follows.
\begin{equation}
\mathcal{L}(\mathbf{x}_L(t))=\hat{e}_{L,L}^\top C_L\mathbf{x}_L(t)
\label{eq:AR_part_Matrix_form}
\end{equation}
One can observe that to each model of the form \eqref{eq:AR-part}, there corresponds a matrix of the form \eqref{eq:C_matrix}. From here on, a matrix that satisfies the previous conditions will be called the matrix form of a linear model $\mathcal{L}$ of the form \eqref{eq:AR-part}.

Given $\delta>0$, and two matrices $A\in \mathbb{C}^{m\times n}$ and $Y\in\mathbb{C}^{m\times p}$, let us write $AX\approx_\delta Y$ to represent the problem of finding $X\in \mathbb{C}^{n\times p}$, $\alpha,\beta\geq 0$ and an orthogonal projector $Q$ such that $\|AX-Y\|_F\leq \alpha\delta+\beta\|(I_m-Q)Y\|_F$. The matrix $X$ will be called a solution to the problem $AX\approx_\delta Y$. 

\begin{thm}\label{thm:thm-1}
Given $\delta>0$, two integers $L,M>0$, a sample $\Sigma_{N}=\{x_t\}_{t=1}^N$ from an AEP signal $\Sigma=\{x_t\}_{t\geq 1}$ with $N>L$, and a matrix $A\in \mathbb{C}^{L\times M}$. If $r=\mathrm{rk}_\delta(\mathcal{H}_L(\Sigma_{N}))>0$, then there is a sparse matrix $\hat{A}\in \mathbb{C}^{L\times M}$ with at most $Mr$ nonzero entries such that $\mathcal{H}_L(\Sigma_N)^\top \hat{A} \approx_\delta \mathcal{H}_L(\Sigma_N)^\top A$.
\end{thm}
\begin{pf}
Since $\mathrm{rk}_\delta(\mathcal{H}_L(\Sigma_{N})^\top)=\mathrm{rk}_\delta(\mathcal{H}_L(\Sigma_{N}))>0$ by \citep[Lemma 3.2]{DBLP:journals/corr/abs-2105-07522}. This result is a consequence of the application of \citep[Theorem 3.6]{DBLP:journals/corr/abs-2105-07522} to the problem $\mathcal{H}_L(\Sigma_N)^\top \hat{A} \approx_\delta \mathcal{H}_L(\Sigma_N)^\top A$.
\end{pf}

Given an AEP signal $\Sigma=\{x_t\}_{t\geq 1}$ with approximate period $T$ and a signal model $\mathcal{S}$ for $\Sigma$ of the form \eqref{eq:ARGRU_model}, if for the linear component $\mathcal{L}$ of $\mathcal{S}$ the corresponding residuals $r_t=|x_{t+1}-\mathcal{L}(\mathbf{x}_L(t))|$ are small, then the significative contribution of $\mathcal{L}$ to the modeling process of $\Sigma$, will be beneficial for interpretability purposes. 

Let us consider an AEP signal $\Sigma$ with approximate period $T$ and approximately periodic tail $\tilde{\Sigma}=\{x_t\}_{t\geq S}$. Given some suitable integer lag value $L>0$ such that there is an approximate linear model $\mathcal{L}$ for $\Sigma$ of the form \eqref{eq:AR-part} with small residuals $r_t=|x_{t+1}-\mathcal{L}(\mathbf{x}_L(t))|$. Let $C_L$ be the matrix form of $\mathcal{L}$ and let $s=S+L-1$. By applying a Krylov subspace approach along the lines presented in \cite[\S 6.1]{doi:10.1137/1.9781611970739}, and as a consequence of \citep[Theorem 4.3.]{DBLP:journals/corr/abs-2105-07522}, one can find a matrix $W_k\in \mathbb{C}^{L\times k}$ whose columns form an orthonormal basis of the subspace $\mathcal{K}_T\subset \mathbb{C}^L$ determined by the expression
\[
\mathcal{K}_{T}=\mathrm{span} \: (\{\mathbf{x}_L(s),C_L \mathbf{x}_L(s),C_L^2 \mathbf{x}_L(s),\ldots,C_L^{T-1}\mathbf{x}_L(s)\}),
\]
such that each $z\in \sigma(W_k^\ast C_LW_k)$ satisfies the relation $|z^T-1|\leq \varepsilon$, for some $\varepsilon>0$. The matrix $W_k^\ast C_LW_k$ will be called an approximately periodic ({\bf AP}) $\Sigma$-section of $C_L$ and will be denoted by $C_L|_{\Sigma}^{AP}$. 

\begin{rem}\label{rem:rem-2}
Based on the previous considerations, when a given AEP signal $\Sigma=\{x_t\}_{t\geq 1}$ with approximate period $T$ is {\em well explained} by the linear component of a semilinear model $\mathcal{S}$ of the form \eqref{eq:ARGRU_model}, that is, when the corresponding residuals are relatively small, the matrix $C_L|_{\Sigma}^{AP}$ corresponding to the model should {\em mimic} the approximate periodicity of the approximately periodic tail $\{x_t\}_{t\geq S}$ of $\Sigma$, in the sense that the number $\|(C_L|_{\Sigma}^{AP})^T-I_k\|$ should be relatively small for some suitable matrix norm $\|\cdot\|$ (in the sense of \cite[\S 1.5]{doi:10.1137/1.9781611970739}). Consequently, the points in $\sigma((C_L|_{\Sigma}^{AP})^T)$ should cluster around $1$.
\end{rem}

\section{Algorithms}\label{sec:algorithm}

The ideas in section \S\ref{sec:semilinear-approximation} can be tanslated into a prototypical algorithm represented by algorithm \ref{alg:main_AutoRegressor_alg_1}, that relies on \cite[Algorithm 1]{DBLP:journals/corr/abs-2105-07522}, Theorem \ref{thm:thm-1} and \cite[Theorem 4.3.]{DBLP:journals/corr/abs-2105-07522}.

\begin{algorithm2e}
\caption{{\bf SpARSModel}: Semilinear Sparse Model Parameters Computation}
\label{alg:main_AutoRegressor_alg_1}
\SetAlgoLined
 \KwData{$\Sigma_{N}=\{x_{t}\}_{t=1}^{N}\subset \mathbb{C}$}
  \KwResult{$\mathbf{c},\mathfrak{G}_1,\ldots,\mathfrak{G}_j,\mathbf{w}_M=\mathbf{SpARSModel}(\Sigma_{N})$}
\begin{itemize}
\item [0:] Estimate the lag value $L$ using auto-correlation function based methods\;
\item [1:] Approximately solve \eqref{eq:sp_autoregressor_eq} for $\mathbf{c}$ using the data in $\Sigma_{N}$ and applying \cite[Algorithm 1]{DBLP:journals/corr/abs-2105-07522}\;
\item[2:] Fit the model blocks $\mathfrak{G}_j$ of \eqref{eq:Model_diagram} using the data in $\Sigma_{N}$.
\item[3:] For the GRU layers of each $\mathfrak{G}_j$, compute sparse representations of the corresponding input weights $W_{ir},W_{iz},W_{in}$ in \eqref{eq:GRU-cell-description} when appropriate, applying Theorem \ref{thm:thm-1}.
\item[4:] Compute the coefficients $\mathbf{w}_M=(w_1,\ldots,w_{N+1})$ of the mixing block $\mathfrak{M}$ of \eqref{eq:Model_diagram} using the data in $\Sigma_{N}$ and \cite[Algorithm 1]{DBLP:journals/corr/abs-2105-07522}\;
\end{itemize}
\KwRet{$\mathbf{c},\mathfrak{G}_1,\ldots,\mathfrak{G}_j,\mathbf{w}_M$}
\end{algorithm2e}
We can apply algorithm \ref{alg:main_AutoRegressor_alg_1} to compute the model parameters needed for the computation of signal models of the form \eqref{eq:Model_diagram}.

\section{Numerical Experiments}
\label{sec:experiments}

In this section, some computational implementations of the methods reported in this document are presented. The experimental results documented in this section can be replicated using the function {\tt NumericalExperiment}{\tt .py} or the Jupyter notebook {\tt SLSpAARModelsDemo.ipynb}, that are available in \cite{FVides_SPAAR}. The configuration required to replicate the results in this section is available as part of the aforementioned programs. 

From here on, we will refer to the sparse semilinear models proposed in this document as {\bf SpARS} models. The signal approximations computed using the SpARS models presented in this document are compared with the approximations obtained using standard AR models, and the standard AR models are computed using the Python program {\tt Autoreg} included as part of the {\tt statsmodels} module. In this section we will write $nz$ to denote nonzero elements. For the experiments documented in this section two GRU RNN blocks were used, the block $\mathfrak{G}_1$ was computed using TensorFlow 2.6.0 and its input weights were replaced by their corresponding sparse representations, that were computed using \cite[Algorithm 1]{DBLP:journals/corr/abs-2105-07522} along the lines of Theorem \ref{thm:thm-1}, and the block $\mathfrak{G}_2$ was computed using PyTorch 1.9.1+cpu and its input weights were left unchanged.

\subsection{Numerical Experiment 1}\label{sec:experiment_1}

\begin{figure}
\centering
\includegraphics[scale=.6]{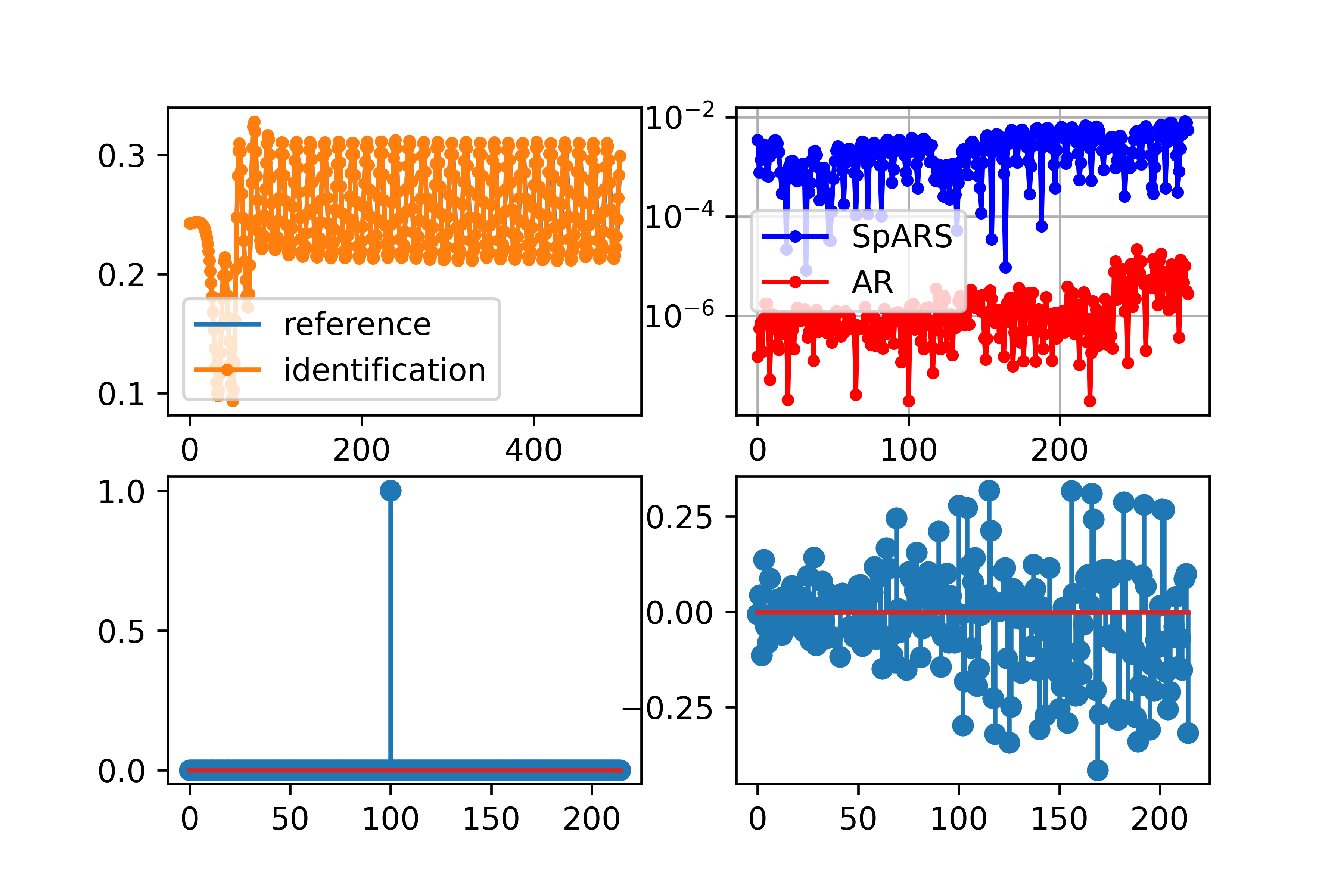}
\caption{Reference and identified signals for SpARS model (top left). Approximation errors (top right). Linear component parameters of the SpARS model with 1 $nz$ (bottom left). Linear component parameters of the standard AR model with 215 $nz$ (bottom right).}
\label{fig:experiment_1_1}
\end{figure}

\begin{figure}
\centering
\includegraphics[scale=.6]{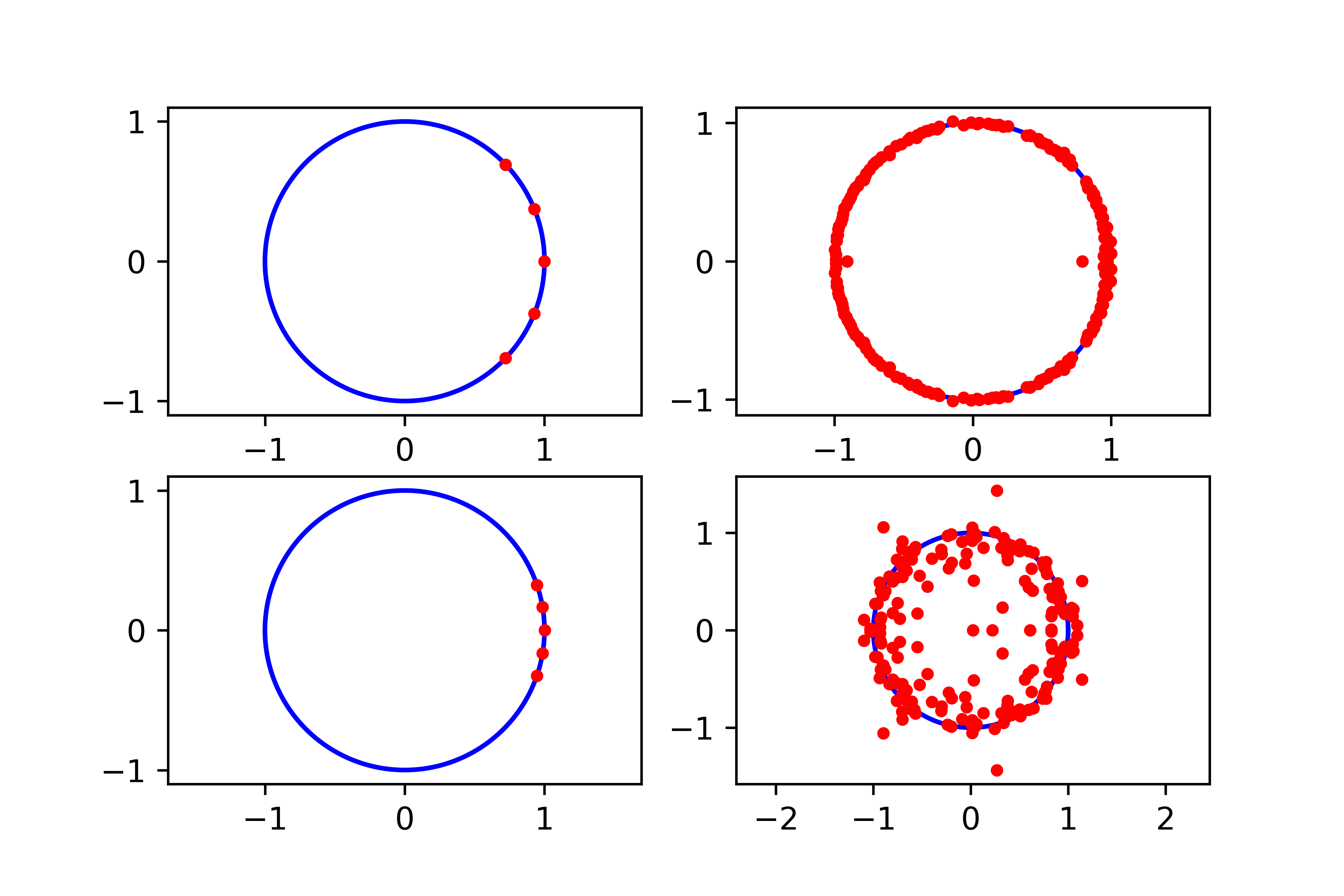}
\caption{$\sigma(C_L|_{\Sigma}^{AP})$ for the linear component of the SpARS model (top left). $\sigma(C_L|_{\Sigma}^{AP})$ for the linear component of the standard AR model (top right). $\sigma((C_L|_{\Sigma}^{AP})^T)$ for the linear component of the SpARS model (bottom left). $\sigma((C_L|_{\Sigma}^{AP})^T)$ for the linear component of the standard AR model (bottom right).}
\label{fig:experiment_1_2}
\end{figure}

In this section, algorithm \ref{alg:main_AutoRegressor_alg_1} is applied to compute a SpARS model for the signal data sample recorded in the csv file {\tt VortexSheddingSignal.csv} in the DataSets folder in \cite{FVides_SPAAR}. The graphic representations of the results produced by the command line
\begin{verbatim}
>>> NumericalExperiment(1)
\end{verbatim}
are shown in figures \ref{fig:experiment_1_1} and \ref{fig:experiment_1_2}, respectively.

In every figure like figure \ref{fig:experiment_1_2}, the red dots represent the points in each considered spectrum, the blue lines represent $\mathbf{S}^1$, and the number $T$ represents the estimated approximate period for each signal considered.

\subsection{Numerical Experiment 2}\label{sec:experiment_2}

\begin{figure}
\centering
\includegraphics[scale=.6]{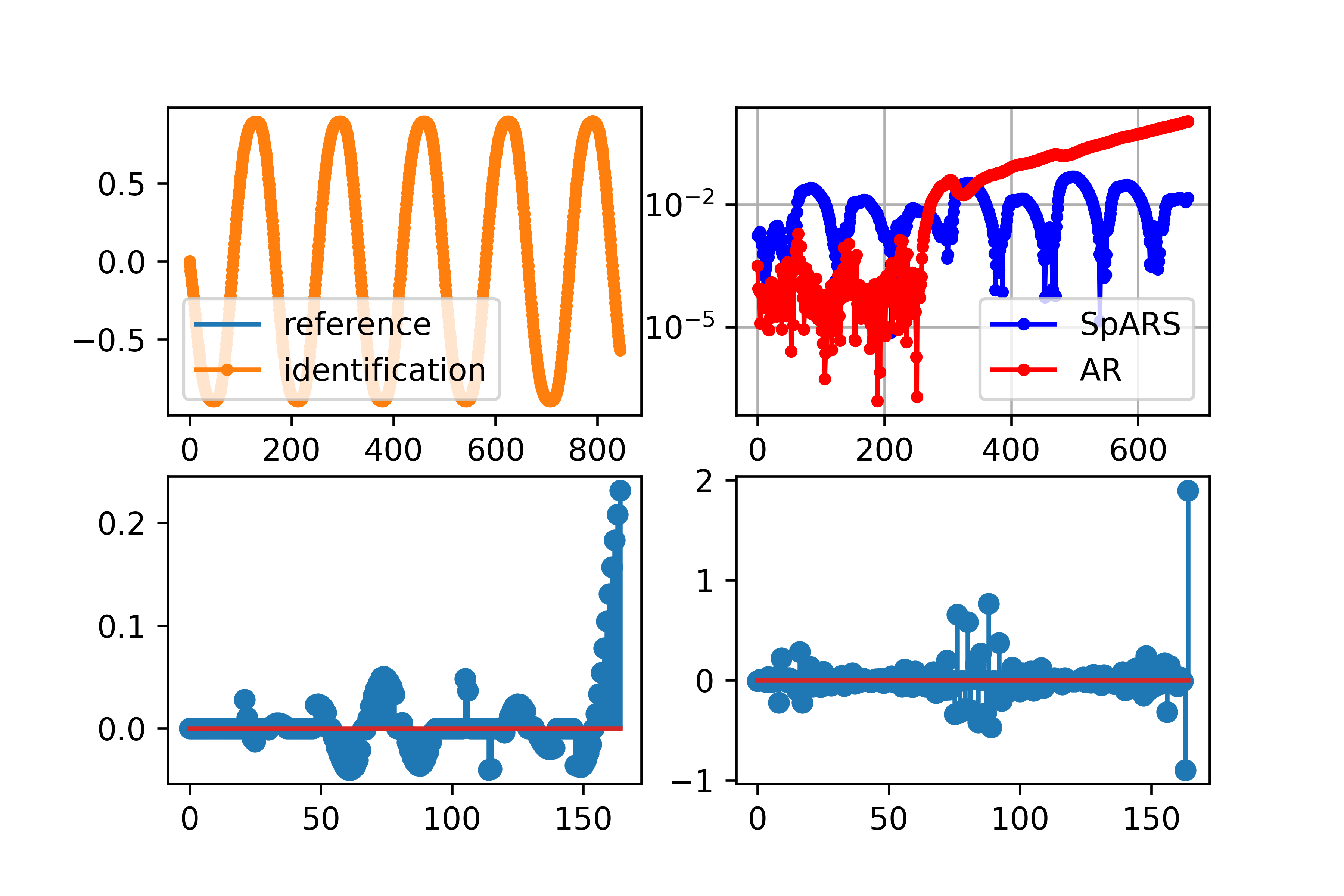}
\caption{Reference and identified signals for SpARS model (top left). Approximation errors (top right). Linear component parameters of the SpARS model with 92 $nz$ (bottom left). Linear component parameters of the standard AR model with 165 $nz$ (bottom right).}
\label{fig:experiment_2_1}
\end{figure}

\begin{figure}
\centering
\includegraphics[scale=.6]{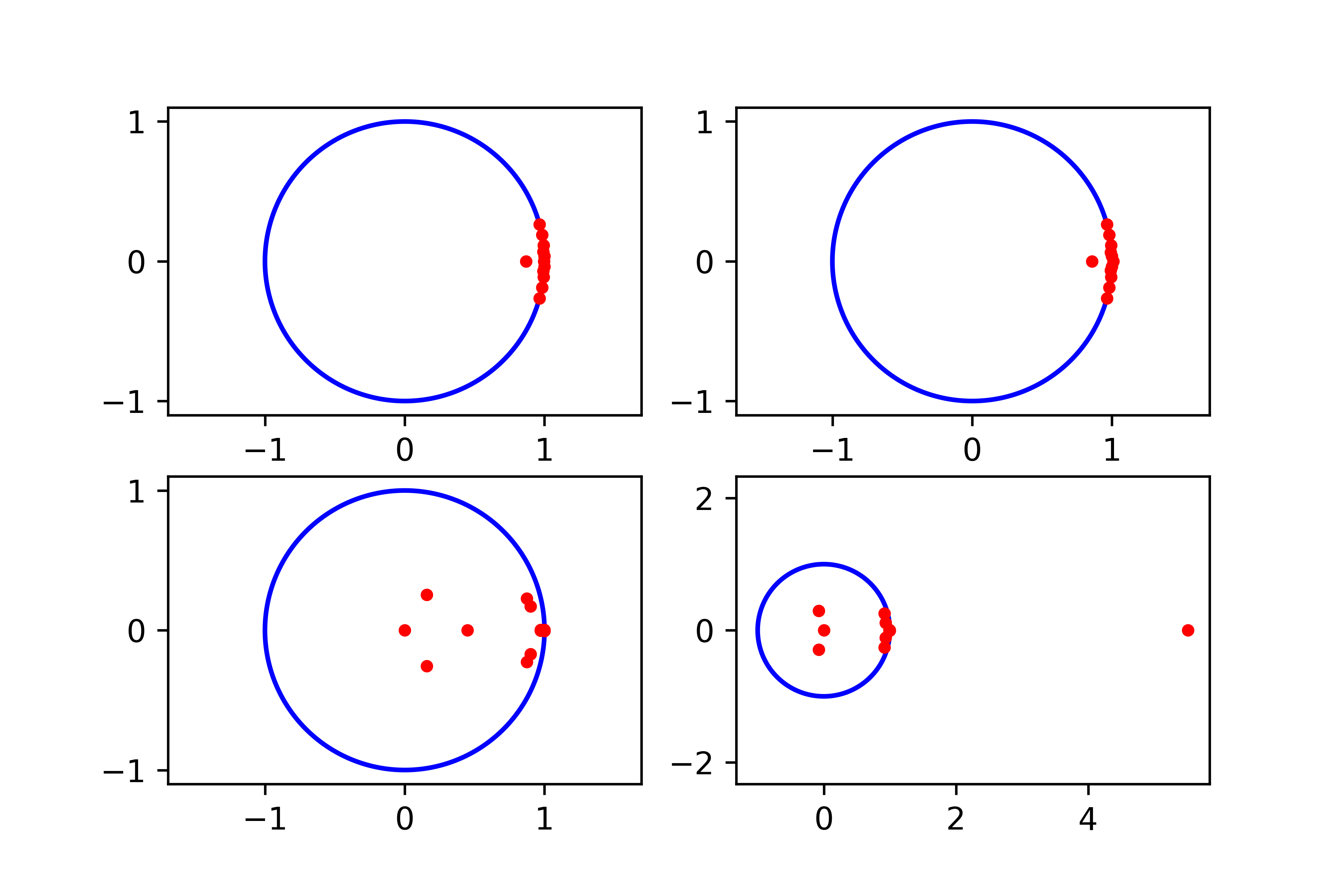}
\caption{$\sigma(C_L|_{\Sigma}^{AP})$ for the linear component of the SpARS model (top left). $\sigma(C_L|_{\Sigma}^{AP})$ for the linear component of the standard AR model (top right). $\sigma((C_L|_{\Sigma}^{AP})^T)$ for the linear component of the SpARS model (bottom left). $\sigma((C_L|_{\Sigma}^{AP})^T)$ for the linear component of the standard AR model (bottom right).}
\label{fig:experiment_2_2}
\end{figure}

In this section, algorithm \ref{alg:main_AutoRegressor_alg_1} is applied to compute a SpARS model for the signal data sample recorded in the csv file {\tt NLOscillatorSignal.csv} in the DataSets folder in \cite{FVides_SPAAR}. The graphic representations of the results produced by the command line
\begin{verbatim}
>>> NumericalExperiment(2)
\end{verbatim}
are shown in figures \ref{fig:experiment_2_1} and \ref{fig:experiment_2_2}, respectively.

\subsection{Numerical Experiment 3}\label{sec:experiment_3}

\begin{figure}
\centering
\includegraphics[scale=.6]{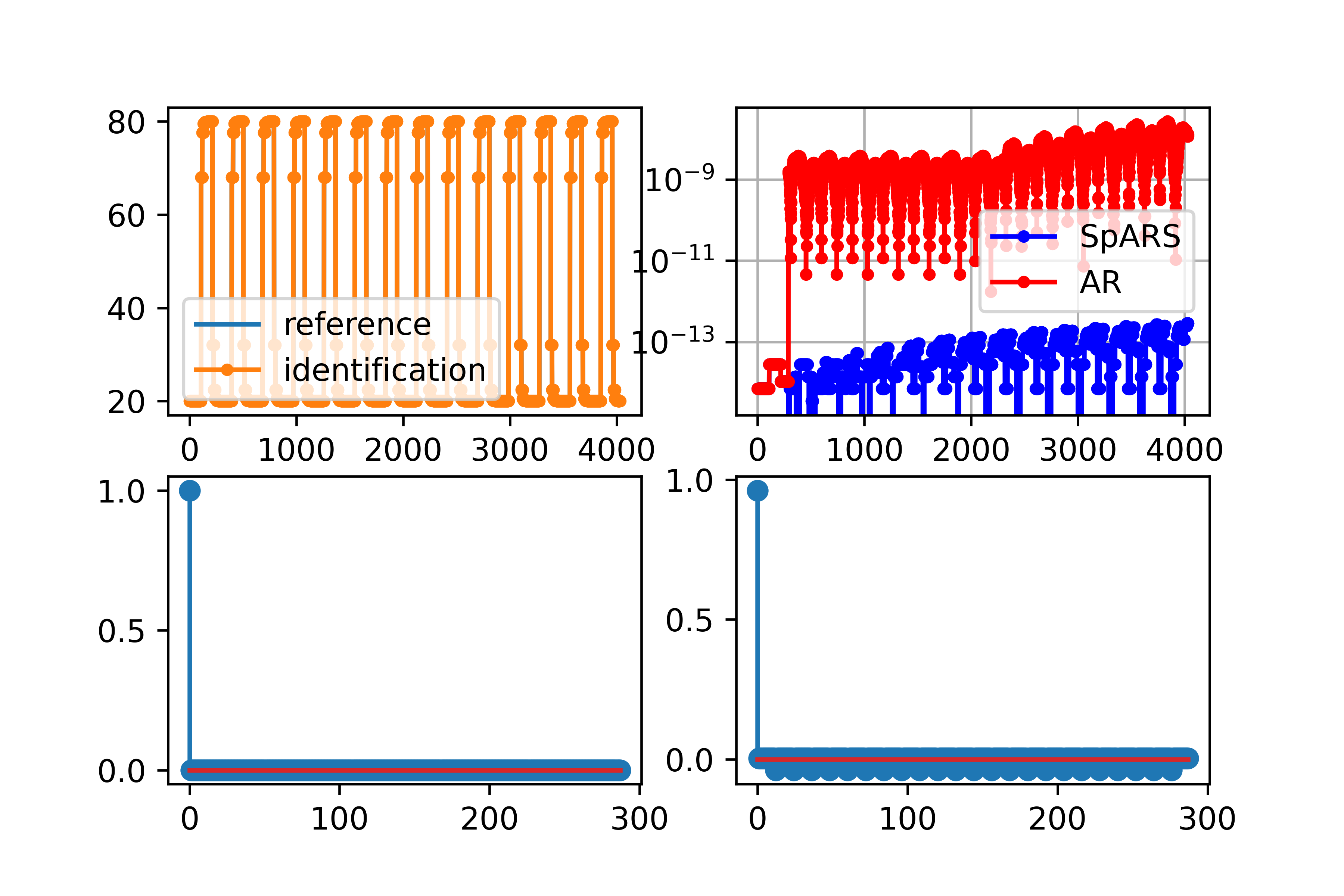}
\caption{Reference and identified signals for SpARS model (top left). Approximation errors (top right). Linear component parameters of the SpARS model with 8 $nz$ (bottom left). Linear component parameters of the standard AR model with 288 $nz$ (bottom right).}
\label{fig:experiment_3_1}
\end{figure}

\begin{figure}
\centering
\includegraphics[scale=.6]{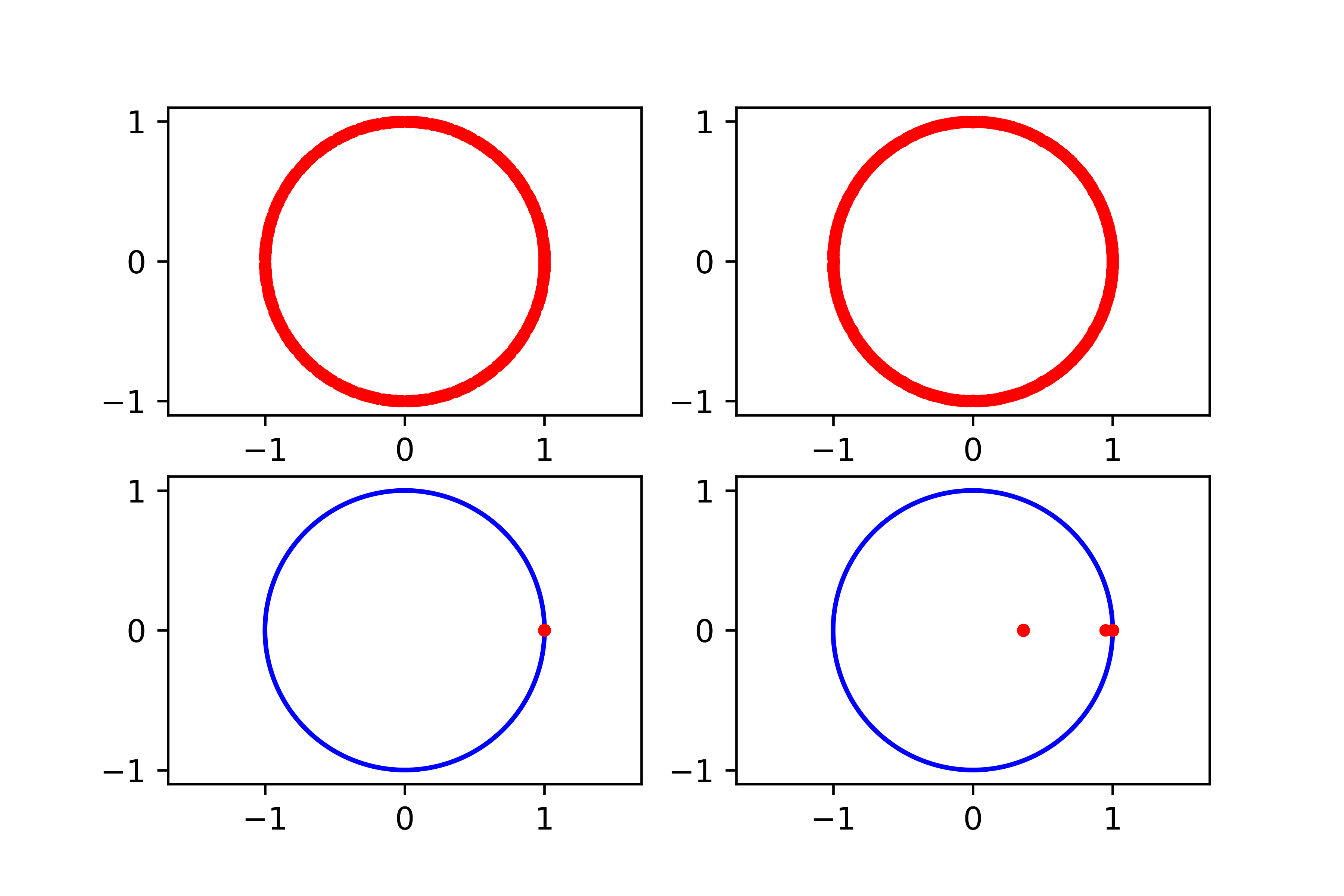}
\caption{$\sigma(C_L|_{\Sigma}^{AP})$ for the linear component of the SpARS model (top left). $\sigma(C_L|_{\Sigma}^{AP})$ for the linear component of the standard AR model (top right). $\sigma((C_L|_{\Sigma}^{AP})^T)$ for the linear component of the SpARS model (bottom left). $\sigma((C_L|_{\Sigma}^{AP})^T)$ for the linear component of the standard AR model (bottom right).}
\label{fig:experiment_3_2}
\end{figure}

\begin{figure}
\centering
\includegraphics[scale=.6]{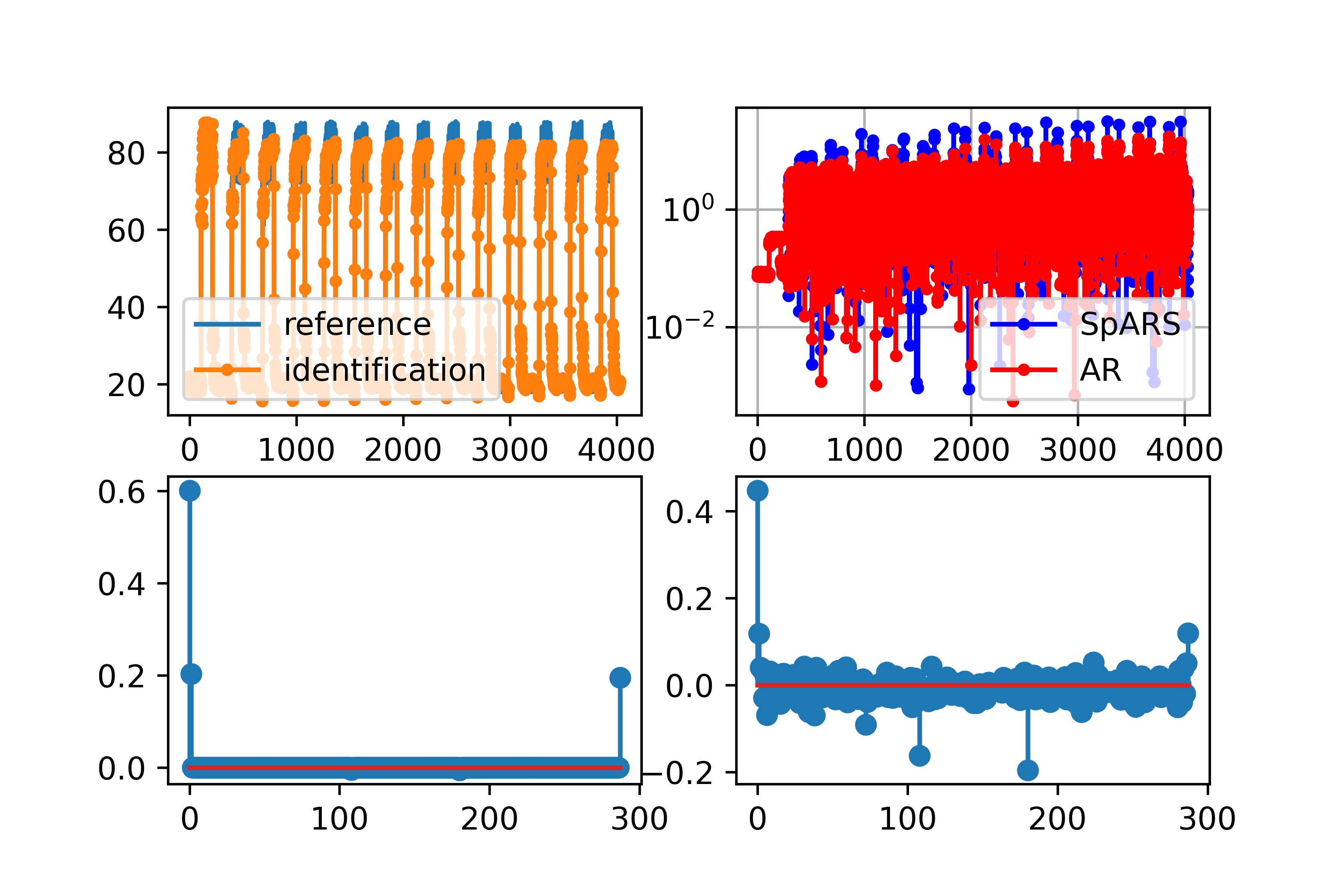}
\caption{Reference and identified signals for SpARS model (top left). Approximation errors (top right). Linear component parameters of the SpARS model with 5 $nz$ (bottom left). Linear component parameters of the standard AR model with 288 $nz$ (bottom right).}
\label{fig:experiment_4_1}
\end{figure}

\begin{figure}
\centering
\includegraphics[scale=.6]{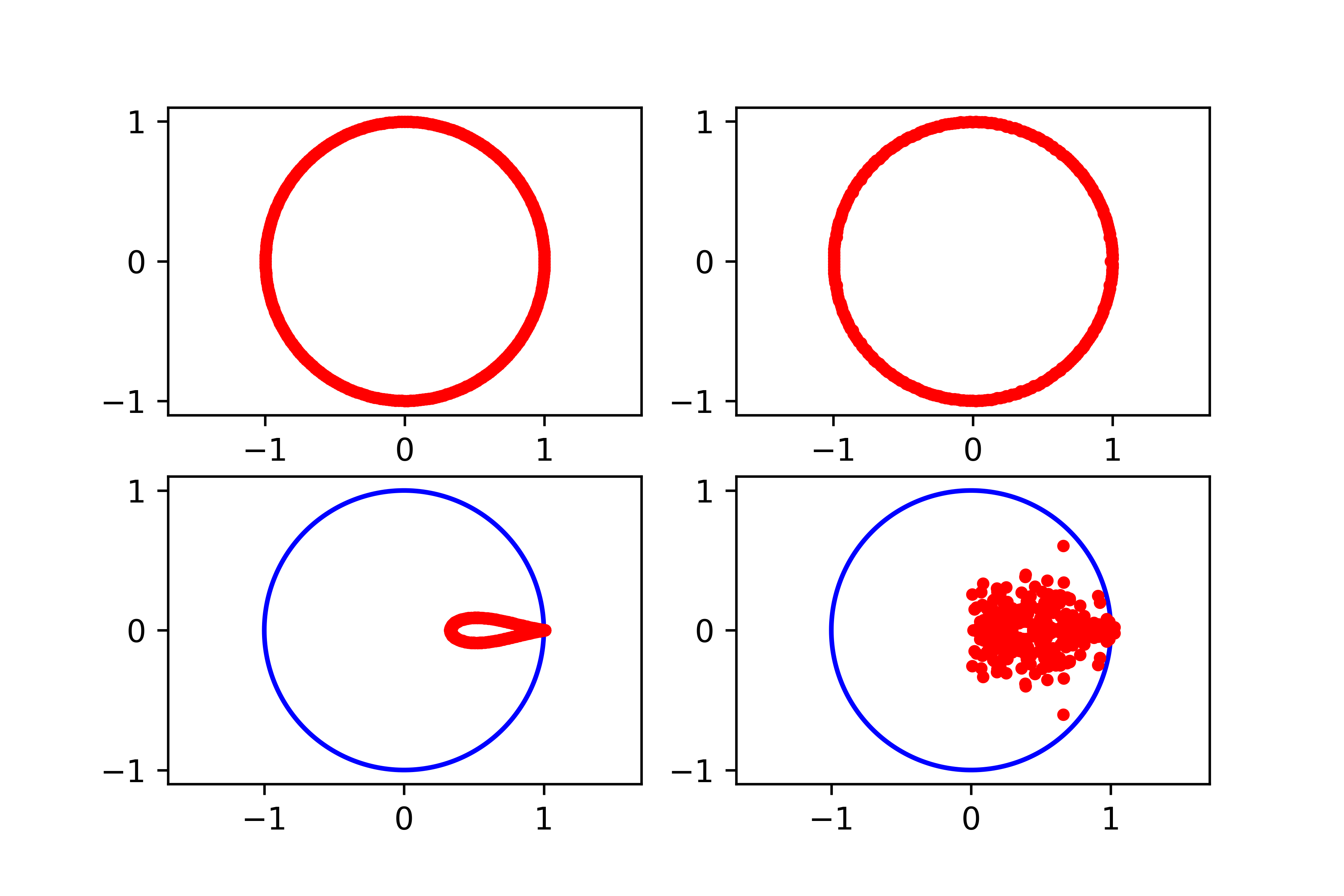}
\caption{$\sigma(C_L|_{\Sigma}^{AP})$ for the linear component of the SpARS model (top left). $\sigma(C_L|_{\Sigma}^{AP})$ for the linear component of the standard AR model (top right). $\sigma((C_L|_{\Sigma}^{AP})^T)$ for the linear component of the SpARS model (bottom left). $\sigma((C_L|_{\Sigma}^{AP})^T)$ for the linear component of the standard AR model (bottom right).}
\label{fig:experiment_4_2}
\end{figure}

In this section, algorithm \ref{alg:main_AutoRegressor_alg_1} is applied to compute a SpARS model for the signal data sample recorded in the csv files: 
\begin{itemize}
\item {\tt art\_daily\_no\_noise.csv}
\item {\tt art\_daily\_small\_noise.csv}
\end{itemize}
that are included as part of the datasets described in \cite{AHMAD2017134}. The graphic representations of the results produced by the command line 
\begin{verbatim}
>>> NumericalExperiment(3.1)
\end{verbatim}
for the periodic signal without noise are shown in figures \ref{fig:experiment_3_1} and \ref{fig:experiment_3_2}, respectively. The graphic representations of the results produced by the command line
\begin{verbatim}
>>> NumericalExperiment(3.2)
\end{verbatim}
for the periodic signal with noise are shown in figures \ref{fig:experiment_4_1} and \ref{fig:experiment_4_2}.

\subsection{Approximation Errors}

The approximation root mean square errors ({\bf RMSE}) are summarized in table \ref{tb:errors}.

\begin{table}[!hb]
\begin{center}
\caption{RMSE}\label{tb:errors}
\begin{tabular}{ccc}
Experiments & SpARS Model & AR Model \\\hline
Experiment 1 & 0.0031265041 & 0.0000035594 \\
Experiment 2 & 0.0162975118 & 0.3100591516 \\
Experiment 3.1 & 0.0000000000 & 0.0000000074 \\
Experiment 3.2 & 3.9894749119 & 4.0939437825 \\ \hline
\end{tabular}
\end{center}
\end{table}

It is appropriate to mention that the root mean square errors can present little fluctuations as one performs several numerical simulations, due primarily to the nature of the neural-network models, as the linear components tend to present very low or no variability from simulation to simulation.

\subsection{Data Availability}

The Python programs that support the findings of this study are openly available in the SPAAR
repository, reference number \cite{FVides_SPAAR}. The time series data used for the experiments 1 and 2 documented in \S\ref{sec:experiment_1} and \S\ref{sec:experiment_2}, respectively, are available as part of \cite{FVides_SPAAR}. The time series data used for experiment 3 in \S\ref{sec:experiment_3} are available as part of  the Numenta Anomaly Benchmark (NAB) described in \cite{AHMAD2017134}.

\section{Conclusion}

Although in some experiments in \S\ref{sec:experiments} the root mean square errors corresponding to the AR and SpARS models are similar, the AP $\Sigma$-sections corresponding to the SpARS models exhibit a better {\em mimetic approximately periodic behavior} in the sense of remark \ref{rem:rem-2}, as it can be visualized in figures \ref{fig:experiment_1_2}, \ref{fig:experiment_2_2}, \ref{fig:experiment_3_2} and \ref{fig:experiment_4_2}. This mimetic behavior is interesting not just from a theoretical point of view, as it provides a criterion for how well the linear component of a given model {\em mimics or captures} the eventual approximate periodic behavior of the signal under study, but also for practical computational reasons, as long term predictions or simulations can be affected when the eigenvalues of the AP $\Sigma$-section of the matrix form corresponding to the linear component of a signal model, lie outside the set $\mathbf{D}^1=\{z\in \mathbb{C}:|z|\leq 1\}$, as one can observe in figures \ref{fig:experiment_2_1} and \ref{fig:experiment_2_2}. Another advantage of the SpARS modeling technology is the reduction of the amount of data needed to fit the corresponding linear component, when compared to a nonsparse linear model, as illustrated in \S\ref{sec:experiment_1}. Although the RMSE for the standard AR model was smaller than the RMSE of the SpARS model, only $50\%$ of the reference data were needed to fit the sparse model, while $90\%$ of the reference data were needed to fit the standard AR model. This difference for the amounts of training data was only necessary for the experiment in \S\ref{sec:experiment_1} due to the relatively small data sample size with respect to the lag value, for the other experiments, $50\%$ of the reference data were used to train both the sparse and the nonsparse models.

\begin{ack}
The numerical simulations documented in this paper were performed with computational resources from the Scientific Computing Innovation Center of UNAH, as part of the researh project PI-063-DICIHT. Some experiments were performed on Google Colab and IBM Quantum Lab as well.
\end{ack}

\bibliography{root}       

\begin{thebibliography}{9}
\providecommand{\natexlab}[1]{#1}
\providecommand{\url}[1]{\texttt{#1}}
\providecommand{\urlprefix}{URL }
\expandafter\ifx\csname urlstyle\endcsname\relax
  \providecommand{\doi}[1]{doi:\discretionary{}{}{}#1}\else
  \providecommand{\doi}{doi:\discretionary{}{}{}\begingroup
  \urlstyle{rm}\Url}\fi

\bibitem[{Ahmad et~al.(2017)Ahmad, Lavin, Purdy, and Agha}]{AHMAD2017134}
Ahmad, S., Lavin, A., Purdy, S., and Agha, Z. (2017).
\newblock Unsupervised real-time anomaly detection for streaming data.
\newblock \emph{Neurocomputing}, 262, 134--147.
\newblock \doi{https://doi.org/10.1016/j.neucom.2017.04.070}.
\newblock Online Real-Time Learning Strategies for Data Streams.

\bibitem[{Cho et~al.(2014)Cho, van Merri{\"e}nboer, Gulcehre, Bahdanau,
  Bougares, Schwenk, and Bengio}]{cho-etal-2014-learning}
Cho, K., van Merri{\"e}nboer, B., Gulcehre, C., Bahdanau, D., Bougares, F.,
  Schwenk, H., and Bengio, Y. (2014).
\newblock Learning phrase representations using {RNN} encoder{--}decoder for
  statistical machine translation.
\newblock In \emph{Proceedings of the 2014 Conference on Empirical Methods in
  Natural Language Processing ({EMNLP})}, 1724--1734. Association for
  Computational Linguistics, Doha, Qatar.
\newblock \doi{10.3115/v1/D14-1179}.
\newblock \urlprefix\url{https://aclanthology.org/D14-1179}.

\bibitem[{Chollet et~al.(2015)}]{chollet2015keras}
Chollet, F. et~al. (2015).
\newblock Keras.
\newblock \url{https://keras.io}.

\bibitem[{Khandelwal et~al.(2015)Khandelwal, Adhikari, and
  Verma}]{KHANDELWAL2015173}
Khandelwal, I., Adhikari, R., and Verma, G. (2015).
\newblock Time series forecasting using hybrid arima and ann models based on
  dwt decomposition.
\newblock \emph{Procedia Computer Science}, 48, 173--179.
\newblock \doi{https://doi.org/10.1016/j.procs.2015.04.167}.
\newblock International Conference on Computer, Communication and Convergence
  (ICCC 2015).

\bibitem[{Paszke et~al.(2019)Paszke, Gross, Massa, Lerer, Bradbury, Chanan,
  Killeen, Lin, Gimelshein, Antiga, Desmaison, Kopf, Yang, DeVito, Raison,
  Tejani, Chilamkurthy, Steiner, Fang, Bai, and Chintala}]{NEURIPS2019_9015}
Paszke, A., Gross, S., Massa, F., Lerer, A., Bradbury, J., Chanan, G., Killeen,
  T., Lin, Z., Gimelshein, N., Antiga, L., Desmaison, A., Kopf, A., Yang, E.,
  DeVito, Z., Raison, M., Tejani, A., Chilamkurthy, S., Steiner, B., Fang, L.,
  Bai, J., and Chintala, S. (2019).
\newblock Pytorch: An imperative style, high-performance deep learning library.
\newblock In H.~Wallach, H.~Larochelle, A.~Beygelzimer, F.~d\textquotesingle
  Alch\'{e}-Buc, E.~Fox, and R.~Garnett (eds.), \emph{Advances in Neural
  Information Processing Systems 32}, 8024--8035. Curran Associates, Inc.

\bibitem[{Saad(2011)}]{doi:10.1137/1.9781611970739}
Saad, Y. (2011).
\newblock \emph{Numerical Methods for Large Eigenvalue Problems}.
\newblock Society for Industrial and Applied Mathematics.
\newblock \doi{10.1137/1.9781611970739}.

\bibitem[{Vides(2021{\natexlab{a}})}]{FVides_SPAAR}
Vides, F. (2021{\natexlab{a}}).
\newblock Spaar: Sparse signal identification python toolset.
\newblock \urlprefix\url{https://github.com/FredyVides/SPAAR}.

\bibitem[{Vides(2021{\natexlab{b}})}]{DBLP:journals/corr/abs-2105-07522}
Vides, F. (2021{\natexlab{b}}).
\newblock Sparse system identification by low-rank approximation.
\newblock \emph{CoRR}, abs/2105.07522.
\newblock \urlprefix\url{https://arxiv.org/abs/2105.07522}.

\bibitem[{Zhang(2003)}]{ZHANG2003159}
Zhang, G. (2003).
\newblock Time series forecasting using a hybrid arima and neural network
  model.
\newblock \emph{Neurocomputing}, 50, 159--175.
\newblock \doi{https://doi.org/10.1016/S0925-2312(01)00702-0}.

\end{thebibliography}

\end{document}